\documentclass[centertags,12pt,a4paper]{amsart}

\usepackage{listings}
\usepackage{color}
\usepackage{latexsym}
\usepackage{amssymb,amsmath}
\usepackage{amscd}

\makeatletter
\@addtoreset{equation}{section}
\makeatother

\newtheorem{theorem}{Theorem}[section]

\newtheorem{definition}[theorem]{Definition}
\newtheorem{lemma}[theorem]{Lemma}

\newtheorem*{claim*}{Claim}
\theoremstyle{definition}

\newtheorem*{question*}{Question}

\def\F{\mathbf{F}}
\def\cB{\mathcal{B}}

\def\cF{\mathcal{F}}

\def\cT{\mathcal{T}}
\def\cC{\mathcal{C}}
\def\cD{\mathcal{D}}
\def\cF{\mathcal{F}}

\def\cI{\mathcal{I}}

\def\cS{\mathcal{S}}

\def\fq{\mathbf{F}_q}

\def\ha{\textstyle\frac{1}{2}}

\author[Faina, Parrettini and Pasticci]{G. Faina, C. Parrettini and F. Pasticci}

\thanks{2000 {\em Math. Subj. Class.}: 51E21}

\title[Hyperfocused arcs in $PG(2,32)$]{Hyperfocused arcs in $PG(2,32)$}

\address{Dipartimento di Matematica e Informatica, Universit\`a degli Studi di Perugia, I-06123
Perugia, Italy}\email{pasticci@dmi.unipg.it}

\begin{document}

\begin{abstract}
In $PG(2,32)$ the following two results are proven by a computer aided search.
\begin{itemize}
\item[(i)] Uniqueness of hyperfocused $12$-arcs, up to projectivities;
\item[(ii)] Non-existence of hyperfocused $14$-arcs.
\end{itemize}
The existence problem for hyperfocused $16$-arcs remains open.
\end{abstract}

\maketitle

     \section{Introduction}\label{s1}

In a $2$-level secret sharing scheme, the secret is shared among persons (also called participants or representatives) qualified in two levels. From the top level, just two persons are necessary and  sufficient to reconstruct the secret (and initiate the action), whereas from the lower level $n$ persons are needed. Furthermore, the secret
can be reconstructed when $n-1$ persons from the lower level are joined by one top level person.

The simplest non-trivial case, $n=3$, was considered by Simmons \cite{SIM2}, see also \cite{SIM1} who introduced the purely geometric concept of a
hyperfocused arc in an attempt to perform such a secret
sharing scheme using finite geometry.
In his setup secrets are represented by points so that lines are used to describe the scheme.

Suppose that the secret is the point $X$ on a given line $s$ in the three-dimensional projective space $PG(3,q)$. The set of shadows (pieces of information given to the participants) is a subset
of points $\{P_1,\ldots,P_m\}$ on a line $\ell$ through $X$ in case of participants of top level, while is a subset $\cS$ of points of a plane $\alpha$ containing $\ell$ in case of participants of lower level. Suppose that $s\neq \ell$ and that $s$ does not lie on $\alpha$. The set $\{X,P_1,\ldots,P_m\}$ is denoted by $\cI$. Furthermore, $\cS$ must be chosen in such a way that no three of its points are collinear, no two of its points collinear with a point from $\cI$, and no point of $\cD$ is on $\ell$. Equivalently, $\cS$ is required to be an arc disjoint from $\ell$ such that no point from $\cI$ lies one secant of $\cS$. The main goal is to construct efficient schemes of this kind, that is, schemes with $|\cS|+|\cI|$ as large as possible, see \cite{SIM2} and \cite{BW}.

It is useful to look at $\alpha$ as the affine plane $AG(2,q)$ whose projective closure is $PG(2,q)$ with respect to the infinite line $\ell$.
Then an efficient Simmon's scheme requires an arc $\cS$ in $AG(2,q)$ with the property
that the direction number $\mu(\cS)$ of $\cS$, that is, the number
of directions determined by secants of $\cS$, is as small as
possible. Obviously, $\mu(\cS)\geq k-1$ where $k$ is the size of
$\cS$. For $q$ odd, Simmons \cite{SIM2} considered the case
$\mu(\cS)=k$ and called such an arc \emph{sharply focused}.
Actually, for $q$ even, $\mu(\cS)=k-1$ can also occur, and if this
is the case then the arc is called \emph{hyperfocused}. From a
result of Bichara and Korchm\'aros  \cite{BK}, hyperfocused arcs
exist if and only if $q$ is even, see also \cite{CH}.

Beutelspacher and Wettl \cite{BW} relaying on a previous work by
Wettl \cite{WET} obtained the following characterization of sharply
focused set. Let $q=p^a$ be odd and let $\cS$ be sharply focused
arc. Then $\cS$ s is a sharply focused arc if and only if
\begin{itemize}
\item $\cS$ is contained in a conic of $AG(2,q)$
\item $\vert \cS \vert =k$ divides $q+1$, $q-1$ or $q$
\item $\cS$ is affinely regular $k$-gon provided that $p^2$ does not divide $k$.
\end{itemize}
For affine regular $k$-gons in $AG(2,q)$ see \cite{KG1974,KG1976,KS}.

Now, the even order $q$ case is discussed.

Wettl \cite{WET} showed the existence of hyperfocused arcs not
contained in a conic. Small hyperfocused arcs  were classified for $k=4,6,8$ by  D.~Drake and K.~Keating \cite{DK}, and independently  by
Cherowitzo and Holder \cite{CH}, and for $k=10$ by Cherowitzo and Holder \cite{CH}, see also Giulietti and Montanucci \cite{GMM}. They exist
(and are projectively unique) for $k=4,8$, and also for $k=6$ when $q$ is a square, and for $k=10$, as well, when $q$ is a cube. All such small
hyperfocused arcs are contained in a hyperoval. Cherowitzo and Holder \cite{CH} classified all
hyperfocused arcs contained in a hyperconic and solved negatively a problem, posed by
Drake and Keating \cite{DK} on the possible sizes of hyperfocused
arcs.

Giulietti and Montanucci \cite{GMM} constructed hyperfocused
translation arcs which are not contained in either a hyperoval or a
subplane. These arcs are complete in the sense that every point in
$AG(2,q)$ is covered by some chord of the arc. Hyperfocused arcs contained in a translation oval have been classified in \cite{KLS}. 

In \cite{GMM} the authors  generalized the above concept of a
hyperfocused arc in the following way. Let $PG(2,q)$ be the
projective plane of order $q$ over $GF(q)$. Let $\cS$ be a $k$-arc,
that is, an arc of size $k$ in $PG(2,q)$. $\cS$ is called a
\emph{generalized hyperfocused arc} when there is a set $\Omega$ of
$k-1$ non-collinear points in $PG(2,q)$ that blocks every secant of
$\cS$. Note that dropping the adjective ``non-collinear'', that is,
assuming $\Omega$ to be contained in a line $\ell$ then $\cS$ is a
hyperfocused arc in the affine plane arising from $PG(2,q)$ by

removing the line $\ell$. An example of an $8$-arc which is
generalized hyperfocused is found in \cite{GMM}. By a theorem of
Aguglia, Korchm\'aros and Siciliano \cite{AKS} no generalized hyperfocused arc
in $PG(2,q)$ is contained in a conic.

Hyperfocused arcs $\cS$ are related with $1$-factorizations of the
complete graph whose vertices are the points of $\cS$. In fact, each
point $P$ from the set $\Omega$ gives rise naturally to a $1$-factor
whose edges are the $k/2$ secants through $P$. Therefore,
hyperfocused arcs can be viewed as $1$-factorizations embedded in
$PG(2,q)$ with $q$ even. This relationship was exploited in \cite{CH} and
\cite{GMM} for the classification of hyperfocused arcs of size
$k\leq 10$.

In this paper we deal with the classification problem of
hyperfocused arcs for $q=32$.
In $PG(2,32)$ the Cherowitzo-Holder classification establishes the existence and projective uniqueness of hyperfocused $k$-arcs
contained in a hyperconic. The values of $k$ that occur are: $k=8,12,16,32,34$. Furthermore, a hyperfocused arc other than a hyperoval or an oval minus a point, has size at most $\ha q$, see \cite{CH}. Therefore, the possibilities for the sizes of hyperfocused arcs in $PG(2,32)$ are $$k=4,8,12,14,16,32,34.$$
{}From the above quoted results, it is still open the existence problem 
 for $k=14$ and also the classification problem for $k=12$ and $k=16$. 
Our results concern hyperfocused arcs of sizes $k=12$
and $k=14$ and are stated in the following theorems.
\begin{theorem}
\label{main12} Up to projectivities, there is a unique hyperfocused arc of size $12$ in $PG(2,32).$
\end{theorem}

\begin{theorem}
\label{main14} No hyperfocused arc of size $14$ exists in
$PG(2,32)$.
\end{theorem}
It should be noted that our approach does not relay on the
relationship between hyperfocused arcs and $1$-factorizations. The
reason is that while the number of $1$-factorizations on the
complete graph on $10$ vertices is ``only'' $396$, the corresponding
number on $12$ (and $14$) vertices is rather huge, being 526,915,620 and  1,132,835,421,602,062,347
respectively, see \cite{DGM} and \cite{KO}. The main tools in the proofs of the above theorems
come from the theory of arcs in $PG(2,q)$ as developed in
\cite[Chapter X]{HIR}, together with an exhaustive computer
search performed in 
$C{++}$.

Our final remark is that the above results can be re-stated in terms of desarguesian nets, see \cite{BFG,D,DK}.

\section{Definitions and Preliminaries}
Let $\fq$ be the finite fields of order $q$. Let $PG(2,q)$
denote the projective plane over $\fq$ equipped with a projective frame of homogeneous coordinates $(X_1, X_2, X_3)$.  As usual, the affine plane $AG(2,q)$ arising from $PG(2,q)$ by removing the line $\ell$ of equations $X_3=0$ is also considered. Conversely, the projective closure of $AG(2,q)$ equipped with non-homogenous coordinates $(X,Y)$ is $PG(2,q)$ with infinite line $\ell$ and $X=X_1/X_3,\,Y=X_2/X_3$.

An \emph{arc} $K$ in $PG(2,q)$
is a set of points no three of which are collinear. If $K$ has size
$k$ it is a \emph{$k$-arc}. A \emph{secant} (or \emph{chord}) of $K$
is a line $\ell$ containing two points from $K$. A \emph{ blocking
set of the secants} of $K$ is a point set $\cB \subset
PG(2,q)\setminus K$ such that $\cB$ meets every secant of
$K$. It is easy to prove that the size of $\cB$ is at least $k-1$.
If this lower bound is attained, $\cB$ is said to have \emph{minimum
size}.

In $PG(2,q)$ with $q=2^e$, a (complete) $(q+2)$-arc is also called a hyperoval. An irreducible conic plus its nucleus is an example of a hyperoval, called a hyperconic. Also, for every $i$ with $\gcd(i,e)=1$, the set $\Omega_i$ of points $P=(a,a^{2^i},1)$ together with $X_\infty=(1,0,0)$ and $Y_\infty=(0,1,0)$ is a hyperoval. Such a hyperoval is called a translation hyperoval, as the affine points of $\Omega_i$ form an orbit of the group of all translations $(X,Y)\to (X+a,Y+a^{2^i})$ with $a\in GF(q)$.

The following theorem which is a useful tool in our investigation on arcs is due to Bichara and Korchm\'aros \cite{BK},
\begin{theorem}
\label{sandro} Let $\Omega$ be a set of $q+2$ points in $PG(2,q)$.
Let $\Phi$ be the set of points $P\in\Omega$ such that every line
containing $P$ intersects $\Omega$ in two points. Then
\begin{itemize}
\item[(i)] If $|\Phi|>2$ then $q$ is even.
\item[(ii)] If $|\Phi|>q/2$ then $\Phi=\Omega$.
\item[(iii)] For each even $q$, there exists an $\Omega$ such that
$|\Phi|=q/2$.
\end{itemize}
\end{theorem}
Hyperfocused arcs, as well as sharply focused arcs, are particular $k$-arcs in $PG(2,q)$ which are defined as follows.  \begin{definition}
In $PG(2,q)$, let $K$ be a $k$-arc with $k\geq 5$ and let $\ell$ be a line disjoint from $K$. The focus set $\cF$ of $K$ on $\ell$ consists of all focuses of $K$ on $\ell$, thats is, $\cF$ is the set all points on $\ell$ covered by the secants of $K$. If $\cF$ has minimum size, that is, $|\cF|=k-1$, the arc $K$ is said to be hyperfocused. If $|\cF|=k$, the arc $K$ is said to be sharply focused.
\end{definition}
Given a hyperfocused arc $K$ on $\ell$, add to $K$ the points of $\ell$ which are not covered by the secants of $K$. The resulting set $\Omega$ has size $q+2$, and $K$ has the property required for $\Phi$ in Theorem \ref{sandro}. As Cherowitzo and Holder \cite{CH} pointed out, (i) and (ii) of Theorem \ref{sandro} have the following corollary.
\begin{theorem}
If $PG(2,q)$ contains a hyperfocused arc then $q$ is even. Furthermore, if $K$ is a hyperfocused arc of size $k$ distinct from a hyperoval, then $k \leq$ $\frac{q}{2}$.
\end{theorem}
For a hyperfocused arc $K$, every focus is the common point of secants but tangents of $K$. Therefore, the following result holds.
\begin{lemma}
\label{lemma1} If $k$-arc in $PG(2,q)$ is hyperfocused, then $k$ is even.
\end{lemma}
The following result is a corollary of \cite[Proposition 2.1]{BFG}.
\begin{theorem}
\label{BFGth}
Every sharply focused $k$-arc in $PG(2,q)$, with $q$ even, is contained in a hyperfocused $(k+1)$-arc with the same focus set $\cF$.
\end{theorem}
Now, we quote some examples of hyperfocused arcs from Holder's thesis \cite{HOL}.
She pointed out that every hyperoval $\Omega$ in $PG(2,q)$, $q$ even, is a hyperfocused arc on any line $\ell$ disjoint from $\Omega$. Furthermore, if  $PG(2,q)$ is viewed as a subplane in $PG(2,q^e)$ then any hyperfocused arc of $PG(2,q)$ on a line $\ell$ also defines a hyperfocused arc of $PG(2,q^e)$ on $\ell$. On the other hand, every 4-arc in $PG(2,q)$ with $q$ even is hyperfocused on the diagonal line of the 4-arc.


Holder's examples were generalized by Giulietti and Montanucci \cite{GMM}. In their work, the first step was to construct hyperfocused arcs contained in a conic using additive sets lying in a conic $PG(2,q)$ with respect to a tangent. They were also able to give a characterization of such hyperfocused arcs. To quote their result, consider the irreducible conic $\cC$ in $PG(2,q)$ in its canonical form $X_1^2=X_2X_3$ with respect to a projective frame of $PG(2,q)$. For a subgroup $T$ of $(\fq,+)$, the subset $\cT$ of $\cC$ consisting of all points $P=(t,t^2,1)$ with $t$ running over $T$ is the \emph{additive set arising from} $\cT$. With this notation, the following result was proven in \cite{GMM}.
\begin{theorem}
Let $q=2^e$ with $e\geq 2$. A set of points $K$, with $\vert K \vert >3$, contained in a conic is hyperfocused on a tangent line $\ell$ to $\cC$ if and only if $K$ is projectively equivalent to the additive set arising from a non-trivial subgroup of $(\fq ,+)$. In particular, the size of $K$ is a power of $2$.
\end{theorem}
After observing that an additive set is an orbit under an elation group of axis $\ell$, Giulietti and Montanucci went on to investigate $k$-arcs (not necessarily contained in a conic) 
with such an elation group. 
Here we report on their results and constructions from \cite{GMM} which provide a large family of hyperfocused arcs not contained in a conic.

Given a pair $A=(a,b)\in \fq \times \fq $ denote $\overline{A}$ the point in $PG(2,q)$ with coordinates $(a,b,1)$ and let $\phi_A$ be the elation
$$\phi_A:(X_1 ,X_2 ,X_3 )\to (X_1 +a_1 X_3, X_2+a_2X_3, X_3 ).$$
Equivalently, in $AG(2,q)$ arising from $PG(2,q)$ by removing $\ell$, $\phi_A$ is a translation.

Let $G$ be an additive subgroup of $\fq \times \fq $ and let $K_G (P)$ be the orbit of an affine point $P$ under the action of the group
$$T_G:= \{ \phi_A \vert A \in G \}.$$
If such an orbit $K_G(P)$ is a $k$-arc then it is a \emph{translation arc}. Clearly, every additive set is a translation arc.
Interestingly, every translation arc is an hyperfocused arc on $\ell$.
This follows from the following result, see \cite{GMM}.
\begin{theorem}
Let $K$ be a translation arc. Then there exists a blocking set of the secants of $K$ of minimum size which is contained in $\ell$.
\end{theorem}
Let $G=\{ (\alpha , \alpha^{2^i} ) \vert \alpha \in H \}$ for a subgroup $H$ of $(\fq,+)$ and $(i,e)=1$ with $e=log_2 q$. Then the above arc $K_G$ is translation-arc contained in the translation hyperoval $\Omega_i$.

Actually, the family of translation arcs are much more large. This follows from the result below, see \cite{GMM}.
\begin{theorem}
\label{GMth}
Let $K_G$ be a translation $k-$arc in $PG(2,q)$ with $q=2^e$ and $e\geq 5$. Assume that there exists a point $\overline{A} \in PG(2,q)$ belonging to no secant of $K_G$. Then the set $K':=K_G \cup \phi_A(K_G)$ is a translation $2k-$arc. If $K_G$ is contained in a hyperoval and $\overline{A}\not\in K_G$, then $K'$ is not contained in any hyperoval.
\end{theorem}

Finally we prove the following two results.
\begin{theorem}\label{ossK}
Let $K$ and $\overline K \subsetneq K$ be two  hyperfocused arcs in $PG(2,q)$ on the same line $\ell$.  Then  $\vert K \vert  \geq 2 \vert \overline{K} \vert $.
\end{theorem}
\begin{proof}
 Each line joining a point of $\mathcal{F}_{K} \backslash \mathcal{F}_{\overline{K}}$ and a point of $\overline{K}$ is a tangent line of $\overline{K}$ and a secant line of  $K$. Then  $\vert K \backslash \overline{K} \vert \geq k$.
\end{proof}
\begin{lemma}\label{ossK1}
Let $K$ and a sharply focused arc in $PG(2,q)$. 
If $q$ is even, then $|K|$ is odd.
\end{lemma}
\begin{proof} Since $K$ is a sharply focused arc, its $k(k-1)/2$ secants cover a subset $\Sigma=\{P_1,\ldots,P_k\}$ of points on a line $\ell$ disjoint from $K$. We show that this is only possible for $|K|$ odd.
Take a point $Q\in K$. Since there are exactly $k-1$ secants of $K$ through $Q$, just one point $P\in \Sigma$ is uncovered by those secants.
Therefore, the line $PQ$ is a tangent to $K$ from $P$; it will called the ``principal tangent'' of $K$ at $Q$. Our aim is to show that the principal tangents of $K$ meet at a point.  Let $ABC$ be any triangle inscribed in $K$. Let $P_1,P_2,P_3$ the common points of $\Sigma$ with the principal tangents to $K$ at the vertices of $ABC$, say
$$P_1=AB\cap \Sigma,\,\, P_2=BC\cap \Sigma,\,\, P_3=CA\cap \Sigma.$$
Segre's Lemma of tangents applied to $ABC$ implies that the product of the slops of the lines $P_1C,\,P_2A$ and $P_3B$ is equal to $1$. From Segre's result \cite[Theorem X]{HIR}, the principal tangents are lines of an envelop of degree $1$, that is, they are lines of a pencil. If $T$ denotes the center of the pencil, all principal tangents of $K$ pass through $T$. Clearly $P\not \in \ell$. Hence, every point in $\Sigma$ is on exactly one principal tangent of $K$.
Therefore, the number of secants of $K$ through a point of $\Sigma$ must be $(k-1)/2$. In particular, $|K|$ is odd.
\end{proof}
We remark that a different proof of Lemma \ref{ossK1} can be given using Lemma \ref{lemma1} and Theorem \ref{BFGth}.


\section{Classification of hyperfocused arcs of size $12$ in PG(2,32)}
We keep up  notation and terminology from the previous sections. Moreover, if $K$ is an arc and $\ell$ is a line disjoint from $K$ such that $K$ is focused on $\ell$  then $\ell$ is called the focus line, the set of all points on $\ell$ covered by the secants of $K$ is called the focus set (of $K$ on $\ell$) and the focus set is denoted by $\mathcal{F}_K$.

Our aim is to show that a every hyperfocused arc of size $12$ in $PG(2,32)$ is necessarily contained in a hyperconic. This is achieved by an exhaustive computer aided search relaying on the following discussions.

In $PG(2,32)$, let $K_0$ be a $8$-arc contained in a hyperfocused arc $K$. By Theorem \ref{ossK}, $K_0$ is not a hyperfocused arc. From Lemma \ref{ossK1}, $K_0$ is neither a sharply focused arc. Therefore, the focus set of $K_0$ on $\ell$ contains at least $9$ points. 
On the other hand, an exhaustive computer aided search shows that no such $8$-arc $K_0$ with $9$ focuses exist. 

Assume that $K_0$ has exactly $10$ focuses, and consider the set $\Phi$ of those focuses of $K_0$ on $\ell$ which are each the common points of at least three secants of $K$. Since $K_0$ has $28$ secants, at least one focus of $K_0$ does not belong to $\Phi$. Actually, if just one focus is off $\Phi$, say $F_0$, then through $F_0$ there are at least
six tangents of $K_0$. Since $K_0$ is contained in the hyperfocused arc $K$, and $F_0$ is a focus of $K$ as well, each tangent $t$ of $K_0$ through $F_0$ is a secant of $K$. So, $t$ provides a point of $K$ outside $K_0$. But this is impossible, as $|K|<14$. This contradiction shows that at least two focuses of $K_0$ on $\ell$ are off $\Phi$.
If $K_0$ has more than $10$ focuses, the above argument still works.

This proves the existence of two focuses of $K_0$, say $F_1$ and $F_2$, such that the number of tangents of $K_0$ through $F_1$ is four, and the same holds for $F_2$.

One two such focuses, $F_1$ and $F_2$, of $K_0$ have been found, there are at most $16$ possible choices for the remaining four points of $K$, namely the points of the grid defined by the two set of tangents through $F_1$ and $F_2$, respectively. So the four points of $K\setminus K_0$ can be obtained by testing all possible $4$-tuples of common points of the tangents of $K$ through either $F_1$ or $F_2$.

The above discussion suggests the following strategy in order to find all hyperfocused arcs of size $12$ in $PG(2,32)$.
The first step is to give a suitable presentation of all projectively non-equivalent $8$-arcs in $PG(2,32)$ which could be contained in a hyperfocused arc of size $12$ in $PG(2,32)$. For this purpose the following lemma is useful.
   \begin{lemma}\label{8-arco_candidato}
Let $K \subseteq PG(2,32)$ be an hyperfocused $12$-arc on the focus line $\ell$. Let  $P_1, P_2, P_3 \in K$ and $\ell_1=\overline{P_2P_3}$, $\ell_2=\overline{P_1P_3}$ and $\ell_3=\overline{P_1P_2}$. Then there exists a projectivity $\phi$ of $PG(2,32)$ acting on the line-set $\{\ell,\ell_1,\ell_2,\ell_3\}$ as follows:
        $$\left\{
        \begin{array}{cccc}
        & \ell  &\longmapsto&  X_3=0\\
        & \ell_1  &\longmapsto&  X_1=0\\
          & \ell_2  &\longmapsto&  X_2=0\\
          & \ell_3  &\longmapsto&  X_1+X_2+X_3=0
        \end{array}
        \right.
        $$
\end{lemma}
\begin{proof}
The line-set $\{\ell,\ell_1,\ell_2,\ell_3 \}$ consists of four lines, no three of which are concurrent. So the assertion follows from the fundamental theorem of projective geometry.
\end{proof}
{}From Lemma \ref{8-arco_candidato}, any hyperfocused $12$-arc is projectively equivalent to an arc satisfying the following conditions:
\begin{equation*}\label{cond12arco}
$$  \begin{itemize}
        \item The focus line has equation $X_3=0$;
        \item $\ell_1 \cap \ell_2 = (0,0,1) \in K$;
        \item $\ell_1 \cap \ell_3 = (0,1,1) \in K$; \;\;\;\;\;\;\;\;\;\;\;\;\;\;\;\;\;\;\;\;\;\;\;\;\;$(\star)$
        \item $\ell_2 \cap \ell_3 = (1,0,1) \in K$;
        \item $(0,1,0)$, $(1,1,0)$ and $(1,0,0)$ are in the focus set.
    \end{itemize}
$$
\end{equation*}
{}
Therefore, an $8$-arc $K_0$ contained in $K$ is given as follows.
$$K_0 = \{ (0,0,1), (0,1,1), (1,0,1), (1,a,1), (c,d,1), (c,e,1), (f,g,1), (f,h,1) \}$$
for some $a,c,d,e,f,g,h \in \F_{32}.$
Now the following theoretic result from Projective Geometry is needed.
\begin{lemma}\label{collineazione}
The map $f: (X_1,X_2,X_3) \mapsto (X_1^2,X_2^2,X_3^2)$, as well as, its powers $$f^i:(X_1,X_2,X_3) \mapsto (X_1^{2^i},X_2^{2^i},X_3^{2^i}))$$ for all $i \in \mathbb{N}$, are collineations of $PG(2,32)$.
\end{lemma}
The map $f^i$ fixes the points $(0,0,1)$, $(1,0,1)$, $(0,1,1)$, $X_{\infty}$, $Y_{\infty}$ and $P_{\infty}$.
Let $w$ be a primitive element of $GF(32)$. The non-zero elements of $GF(32)$ are partitioned into seven orbits under the action of the automorphism group of $GF(32)$. A representative system of such orbits is the set
\begin{equation}
\label{repr} S=\{ 1, \omega ,\omega^3, \omega^5, \omega^7, \omega^{11}, \omega^{15} \}
\end{equation}
By Lemma \ref{collineazione}, we may assume that $a\in S\setminus \{1\}$. So, the number of cases to investigate all possible such $8$-arcs $K_0$ is equal to $7\cdot 32^6$.

The above discussion gives rise a direct algorithm to deal with all such cases, performed in 
$C{++}$.
Up to projectivities, the exhaustive search produced sixty hyperfocused $12$-arcs $K$ in $PG(2,32)$.

Actually, each of such hyperfocused arc of size $12$ turned out to be contained in a hyperconic. To show this,
the following algorithm was used. After constructing the list $L$ of the $12$ points of a such an arc $K$, compute an equation of the conic $C$ passing through the first 5 points in $L$ and check that each of the remaining $7$ points in $L$ satisfy the same equation or coincides with the nucleus of $C$.

\cite[Theorem 4.5]{CH} establishes the projective uniqueness of hyperfocused arc of size $12$ contained in a hyperoval. 
The proof also provides
a construction of a such an arc using a cyclic subgroup of $PGL(3,q)$ of order $11$. Here we give an alternative construction based on affinely regular polygons, in the spirit of the papers \cite{BW,WET}. A survey of results on affinely regular polygons are found in \cite{KG1976} and in \cite{KS}.

For this purpose, consider $PG(2,32)$ as the projective closure of the affine plane $AG(2,32)$. Let $\cC$ be the conic of equation $X^2+XY+Y^2=1$. Since the polynomial $f(x)=x^2+x+1$ is irreducible over $GF(32)$, the conic $\cC$ is an ellipsis. Its nucleus is the origin $O=(0,0)$. To every point $P=(x,y)$ of $AG(2,q)$ there is associated an element $x+iy$ in the quadratic extension $GF(q^2)=GF(q)(i)$ arising from $f(x)$ such that $f(i)=0$. Note that if $z=x+iy$ then $z^q=(x+y)+iy$. Therefore, the points in $\cC$ are those elements $z=x+iy$ which satisfy the equation $z^{q+1}=1$. In terms of norms,  $P(x,y)\in \cC$ if and only if $\parallel z \parallel=1$. Now, the elements of $GF(q^2)$ with norm $1$ form a cyclic subgroup $Z_{33}$ of order $33$. Take an element $w$ of order $11$ from $Z_{33}$, and consider the points $P_0,P_1,\ldots,P_j,\ldots, P_{10}$ that correspond to the powers $w^j$ of $w$. 

The polygon $P_0P_1\ldots P_{10}$ is an affinely regular eleven-gon in $AG(2,32)$. Let $A_0A_1\ldots A_{10}$ be the regular eleven-gon in the Euclidean plane. The natural one-to-one correspondence $P_j\leftrightarrow A_j$ preserves the parallelisms between chords (including sides). Therefore, the $11$-arc $\{P_0,P_1\ldots P_{10}\}$ have exactly ten focuses on the line at infinity. Also, each focus is the point an infinity of a line through the origin and a point of the arc. Adding $O$ to the $11$-arc $\{P_0,P_1\ldots P_{10}\}$ provides the required $12$-arc hyperfocused on the line at infinity.

\section{Non-existence of hyperfocused 14-arcs in PG(2,32)}
We discuss how to adapt the  method explained in the preceding section for the classification hyperfocused $14$-arcs in $PG(2,32)$. Again, the key idea is to investigate $8$-arcs contained in a hyperfocused $14$-arc; this time we have to consider $8$-arcs with at most $13$ (and at least 9) focuses.


With same setup as in the preceding section, we start considering such an $8$-arc $K_0$, called a candidate $8$-arc, and we try to extend it adding more 6 points. 
Let
$$K_0 = \{(0,0,1),(0,1,1),(1,0,1),(1,a,1),(c,d,1),(c,e,1),(f,g,1),(f,h,1)\}$$
be a candidate $8$-arc contained in a hyperfocused $14$-arc with at most $13$ focuses.
{}From Lemma \ref{collineazione}, we limit ourselves to take and investigate $a$ for $a\in S$, see  (\ref{repr}).

For each candidate $8$-arc, add six points so that the 63 new secant lines add new points to the focus set, so that it has 13 elements. The algorithm is the following:
\begin{itemize}

\item add 6 points such that we have 63 new secant lines ($6 \times 8$ are secant lines such that each of them contains one point among the six new points and one among the 8 old points) give new focuses that, added to the previous, we have $|\cF|=13$;
\item for each $8$-arc  find all focuses $Q_1, \ldots ,Q_s $ such that each of them belongs to exactly 6 tangent lines;
\item select 2 focuses $Q_i$ and $Q_j$ in $\{ Q_1, \ldots ,Q_s  \}$.
The lines joining $Q_i$ and $Q_j$ to the points of $K_0$ meet in 36 points $P_1, \ldots , P_{36}$.
Among them choose a point that belongs to one line to each focus
;
\item take all 6-tuples $(P_{j,1}, P_{j,2}, P_{j,3}, P_{j,4}, P_{j,5}, P_{j,6})$ such that each point belongs to one line for each focus
;
\item test if adding each of  these 6-uples to the arc we obtain an arc again and if it and defines 13 focuses i.e. if it is an hyperfocused 14-arc.
\end{itemize}
An exhaustive computer aided search based on the above algorithm did not produce any hyperfocused $14$-arc containing $K_0$. Therefore we have the following result.

\begin{theorem}
In $PG(2,32)$ there are no hyperfocused $14$-arc.
\end{theorem}

\lstset{language=c++, breaklines=true, breakatwhitespace=false, breakautoindent = true, basicstyle=\tiny}
\lstset{columns=flexible}
 \lstset{tabsize=3}


\begin{lstlisting}
#include <iostream>
#include <fstream>
#include <GaloisField.h>
#include <GaloisFieldElement.h>
#include <list>
#include <algorithm>
#include <string>

using namespace galois;
using namespace std;

struct Insieme_E
{
bool valido;
GaloisFieldElement fuoco;
GaloisFieldElement valori[7];
};

struct Candidato
{
bool valido;
GaloisFieldElement x;
GaloisFieldElement y;
};

bool conta_validi(Candidato lista[])
{
int temp=0;
for (int i=0; i<49; i++)
{
if (lista[i].valido)
temp++;
}
return (temp > 5);
}

void print(string A, string B, string C, unsigned int i, double d)
{
ofstream out("Elenco_archi.txt", ios::out|ios::app);
if (!out)
{
cout << "\n\nImpossibile accedere al file di output !\nStatistiche non inserite.\n\n";
exit(0);
}
out << A << i << B << d << C;
out.close();
}

void print(GaloisFieldElement zero1, GaloisFieldElement uno1, GaloisFieldElement a1, GaloisFieldElement c1, GaloisFieldElement d1, GaloisFieldElement e1, GaloisFieldElement f1, GaloisFieldElement g1, GaloisFieldElement h1, GaloisFieldElement P1x, GaloisFieldElement P1y, GaloisFieldElement P2x, GaloisFieldElement P2y, GaloisFieldElement P3x, GaloisFieldElement P3y, GaloisFieldElement P4x, GaloisFieldElement P4y, GaloisFieldElement P5x, GaloisFieldElement P5y, GaloisFieldElement P6x, GaloisFieldElement P6y, unsigned int *tot)
{
ofstream out("Elenco_archi.txt", ios::out|ios::app);
if (!out)
{
cout << "\n\nImpossibile creare il file di output !\nOperazione annullata.\n\n";
exit(0);
}
out << "{ (" << zero1.index() << "," << zero1.index() << ",1) (" << zero1.index() << "," << uno1.index() << ",1) (" << uno1.index() << "," << zero1.index() << ",1) (" << uno1.index() << "," << a1.index() << ",1) (" << c1.index() << "," << d1.index() << ",1) (" << c1.index() << "," << e1.index() << ",1) (" << f1.index() << "," << g1.index() << ",1) (" << f1.index() << "," << h1.index() << ",1) (" << P1x.index() << "," << P1y.index() << ",1) (" << P2x.index() << "," << P2y.index() << ",1) (" << P3x.index() << "," << P3y.index() << ",1) (" << P4x.index() << "," << P4y.index() << ",1) (" << P5x.index() << "," << P5y.index() << ",1) (" << P6x.index() << "," << P6y.index() << ",1) }\n";
(*tot)++;
out.close();
}

bool found(GaloisFieldElement item, GaloisFieldElement insieme[])
{
for (int i=0;i<7;i++)
{
if (insieme[i] == item)
return true;
}
return false;
}

bool test_secanti(GaloisFieldElement zero1, GaloisFieldElement uno1, GaloisFieldElement a1, GaloisFieldElement c1, GaloisFieldElement d1, GaloisFieldElement e1, GaloisFieldElement f1, GaloisFieldElement g1, GaloisFieldElement h1, GaloisFieldElement L[], int *tot)
{
list<GaloisFieldElement> P_impropri;
P_impropri.push_back(zero1);
P_impropri.push_back(uno1);
P_impropri.push_back(a1);
P_impropri.push_back(d1/c1);
P_impropri.push_back(e1/c1);
P_impropri.push_back(g1/f1);
P_impropri.push_back(h1/f1);
P_impropri.push_back(a1+uno1);
P_impropri.push_back((d1+uno1)/c1);
P_impropri.push_back((e1+uno1)/c1);
P_impropri.push_back((g1+uno1)/f1);
P_impropri.push_back((h1+uno1)/f1);
P_impropri.push_back((d1)/(c1+uno1));
P_impropri.push_back((e1)/(c1+uno1));
P_impropri.push_back((g1)/(f1+uno1));
P_impropri.push_back((h1)/(f1+uno1));
P_impropri.push_back((a1+d1)/(c1+uno1));
P_impropri.push_back((a1+e1)/(c1+uno1));
P_impropri.push_back((a1+g1)/(f1+uno1));
P_impropri.push_back((a1+h1)/(f1+uno1));
P_impropri.push_back((d1+g1)/(c1+f1));
P_impropri.push_back((d1+h1)/(c1+f1));
P_impropri.push_back((e1+g1)/(c1+f1));
P_impropri.push_back((e1+h1)/(c1+f1));
P_impropri.sort();
P_impropri.unique();
int n=P_impropri.size()+1;
(*tot)=n;
if ((n>10) && (n<14))
{
list<GaloisFieldElement>::iterator p=P_impropri.begin();
int i=0;
while (p!=P_impropri.end())
{
L[i]=*(p);
i++;
p++;
}
return true;
}
return false;
}

bool test_iperfocalizzato(GaloisFieldElement zero1, GaloisFieldElement uno1, GaloisFieldElement a1, GaloisFieldElement c1, GaloisFieldElement d1, GaloisFieldElement e1, GaloisFieldElement f1, GaloisFieldElement g1, GaloisFieldElement h1, GaloisFieldElement P1x, GaloisFieldElement P1y, GaloisFieldElement P2x, GaloisFieldElement P2y, GaloisFieldElement P3x, GaloisFieldElement P3y, GaloisFieldElement P4x, GaloisFieldElement P4y, GaloisFieldElement P5x, GaloisFieldElement P5y, GaloisFieldElement P6x, GaloisFieldElement P6y, GaloisFieldElement elenco[], int fuochi)
{
list<GaloisFieldElement> Temp;
Temp.push_back(P1y/P1x);
Temp.push_back(P2y/P2x);
Temp.push_back(P3y/P3x);
Temp.push_back(P4y/P4x);
Temp.push_back(P5y/P5x);
Temp.push_back(P6y/P6x);
Temp.push_back((P1y + uno1)/P1x);
Temp.push_back((P2y + uno1)/P2x);
Temp.push_back((P3y + uno1)/P3x);
Temp.push_back((P4y + uno1)/P4x);
Temp.push_back((P5y + uno1)/P5x);
Temp.push_back((P6y + uno1)/P6x);
Temp.push_back(P1y/(P1x + uno1));
Temp.push_back(P2y/(P2x + uno1));
Temp.push_back(P3y/(P3x + uno1));
Temp.push_back(P4y/(P4x + uno1));
Temp.push_back(P5y/(P5x + uno1));
Temp.push_back(P6y/(P6x + uno1));
Temp.push_back((a1 + P1y)/(P1x + uno1));
Temp.push_back((a1 + P2y)/(P2x + uno1));
Temp.push_back((a1 + P3y)/(P3x + uno1));
Temp.push_back((a1 + P4y)/(P4x + uno1));
Temp.push_back((a1 + P5y)/(P5x + uno1));
Temp.push_back((a1 + P6y)/(P6x + uno1));
Temp.push_back((d1 + P1y)/(c1 + P1x));
Temp.push_back((d1 + P2y)/(c1 + P2x));
Temp.push_back((d1 + P3y)/(c1 + P3x));
Temp.push_back((d1 + P4y)/(c1 + P4x));
Temp.push_back((d1 + P5y)/(c1 + P5x));
Temp.push_back((d1 + P6y)/(c1 + P6x));
Temp.push_back((e1 + P1y)/(c1 + P1x));
Temp.push_back((e1 + P2y)/(c1 + P2x));
Temp.push_back((e1 + P3y)/(c1 + P3x));
Temp.push_back((e1 + P4y)/(c1 + P4x));
Temp.push_back((e1 + P5y)/(c1 + P5x));
Temp.push_back((e1 + P6y)/(c1 + P6x));
Temp.push_back((g1 + P1y)/(f1 + P1x));
Temp.push_back((g1 + P2y)/(f1 + P2x));
Temp.push_back((g1 + P3y)/(f1 + P3x));
Temp.push_back((g1 + P4y)/(f1 + P4x));
Temp.push_back((g1 + P5y)/(f1 + P5x));
Temp.push_back((g1 + P6y)/(f1 + P6x));
Temp.push_back((h1 + P1y)/(f1 + P1x));
Temp.push_back((h1 + P2y)/(f1 + P2x));
Temp.push_back((h1 + P3y)/(f1 + P3x));
Temp.push_back((h1 + P4y)/(f1 + P4x));
Temp.push_back((h1 + P5y)/(f1 + P5x));
Temp.push_back((h1 + P6y)/(f1 + P6x));
Temp.push_back((P1y + P2y)/(P1x + P2x));
Temp.push_back((P1y + P3y)/(P1x + P3x));
Temp.push_back((P1y + P4y)/(P1x + P4x));
Temp.push_back((P1y + P5y)/(P1x + P5x));
Temp.push_back((P1y + P6y)/(P1x + P6x));
Temp.push_back((P2y + P3y)/(P2x + P3x));
Temp.push_back((P2y + P4y)/(P2x + P4x));
Temp.push_back((P2y + P5y)/(P2x + P5x));
Temp.push_back((P2y + P6y)/(P2x + P6x));
Temp.push_back((P3y + P4y)/(P3x + P4x));
Temp.push_back((P3y + P5y)/(P3x + P5x));
Temp.push_back((P3y + P6y)/(P3x + P6x));
Temp.push_back((P4y + P5y)/(P4x + P5x));
Temp.push_back((P4y + P6y)/(P4x + P6x));
Temp.push_back((P5y + P6y)/(P5x + P6x));
list<GaloisFieldElement>::iterator p=Temp.begin();
int nuovi=0;
while (p!= Temp.end())
{
bool found=false;
for (int i=0; i< (fuochi-1);i++)
{
if (*p == elenco[i])
{
found=true;
break;
}
}
if (!found)
nuovi++;
p++;
}
return ((nuovi + fuochi) == 13 );
}

int main()
{
time_t start, end;
time(&start);
double diff;
unsigned int n_archi=0;
unsigned int prim_poly[6]={1,0,1,0,0,1}; //coefficienti del polinomio primitivo, ordinati dal grado minore x^5+x^2+1
GaloisField gf(5,prim_poly);
GaloisFieldElement zero(&gf,0); //il secondo parametro passato è SEMPLICEMENTE un'etichetta ! Da varie prove, però, risulta che i primi due elementi inizializzati, con etichetta zero e uno, corrispondono effettivamente agli omonimi elementi del campo.
GaloisFieldElement uno(&gf,1);
cout << "\nInizio i cicli\n";
for (int i_a=1;i_a<32;i_a++)
{
GaloisFieldElement a(&gf,i_a);
if ((a.index()==0)||(a.index()==1)||(a.index()==3)||(a.index()==5)||(a.index()==7)||(a.index()==11)||(a.index()==15))
{
cout << "a = " << a.index() << "\n";
for (int i_c=2;i_c<32;i_c++)
{
GaloisFieldElement c(&gf,i_c);
cout << "\tc = " << c << "\n";
list<GaloisFieldElement> S1;
if (S1.empty() == false) S1.clear();
S1.push_back(a*c);
S1.push_back(c+uno);
S1.push_back(a*c+c+uno);
S1.sort();
S1.unique();
for (int i_d=1;i_d<32;i_d++)
{
GaloisFieldElement d(&gf,i_d);
int n_d=0;
n_d=count(S1.begin(), S1.end(), d);
if (n_d == 0)
{
for (int i_e=2;i_e<32;i_e++)
{
GaloisFieldElement e(&gf,i_e);
if (e.index() > d.index())
{
int n_e=0;
n_e=count(S1.begin(), S1.end(),e);
if (n_e == 0)
{
for (int i_f=2;i_f<32;i_f++)
{
GaloisFieldElement f(&gf,i_f);
if (f.index() > c.index())
{
list<GaloisFieldElement> S2;
if (S2.empty() == false) S2.clear();
S2.push_back(a*f);
S2.push_back((d*f)/c);
S2.push_back((e*f)/c);
S2.push_back(f+uno);
S2.push_back(a*f+f+uno);
S2.push_back((c+d*f+f)/c);
S2.push_back((c+e*f+f)/c);
S2.push_back((d*f+d)/(c+uno));
S2.push_back((e*f+e)/(c+uno));
S2.push_back((a*c+a*f+d*f+d)/(c+uno));
S2.push_back((a*c+a*f+e*f+e)/(c+uno));
S2.sort();
S2.unique();
for (int i_g=1;i_g<32;i_g++)
{
GaloisFieldElement g(&gf,i_g);
int n_g=0;
n_g=count(S2.begin(), S2.end(), g);
if (n_g == 0)
{
for (int i_h=1;i_h<32;i_h++)
{
GaloisFieldElement h(&gf, i_h);
if (h.index() > g.index())
{
int n_h=0;
n_h=count(S2.begin(), S2.end(), h);
if (n_h == 0)
{
GaloisFieldElement lista[14];
int n_fuochi=0;
if (test_secanti(zero, uno, a, c, d, e, f, g, h, lista, &n_fuochi))
{
bool dimensione_valida=false;
int i_E=0;
Insieme_E E[14];
for (int j=0;j<14;j++)
E[j].valido=false;
for (int i_lista=0;i_lista<(n_fuochi-1);i_lista++)
{
list<GaloisFieldElement> Temp;
list<GaloisFieldElement>::iterator p;
if (Temp.empty() == false)
Temp.clear();
Temp.push_back(zero);
Temp.push_back(uno);
Temp.push_back(lista[i_lista]);
Temp.push_back(a+(lista[i_lista]));
Temp.push_back(d+c * lista[i_lista]);
Temp.push_back(e+c * lista[i_lista]);
Temp.push_back(g+f * lista[i_lista]);
Temp.push_back(h+f * lista[i_lista]);
Temp.sort();
Temp.unique();
p=Temp.begin();
if (Temp.size() == 7)
{
dimensione_valida=true;
E[i_E].valido=true;
E[i_E].fuoco=lista[i_lista];
int j_E=0;
while (p!=Temp.end())
E[i_E].valori[j_E++]=*(p++);
i_E++;
}
}
if (dimensione_valida)
{
Candidato elenco[49];
for (int j_elenco=0;j_elenco<49;j_elenco++)
elenco[j_elenco].valido=false;
int i_elenco=0;
GaloisFieldElement x_temp;
GaloisFieldElement y_temp;
for (int i_E1=0; i_E1<7;i_E1++)
{
for (int i_E2=0;i_E2<7;i_E2++)
{
x_temp=((E[0].valori[i_E1]+E[1].valori[i_E2])/(E[0].fuoco+E[1].fuoco));
y_temp=(E[0].valori[i_E1]+(E[1].fuoco * x_temp));
if (!(((x_temp == zero)&&(y_temp == zero)) || ((x_temp == zero)&&(y_temp == uno)) || ((x_temp == uno)&&(y_temp == zero)) || ((x_temp == uno)&&(y_temp == a)) || ((x_temp == c)&&(y_temp == d)) || ((x_temp == c)&&(y_temp == e)) || ((x_temp == f)&&(y_temp == g)) || ((x_temp == f)&&(y_temp == h))))
{
elenco[i_elenco].x=x_temp;
elenco[i_elenco].y=y_temp;
elenco[i_elenco].valido=true;
i_elenco++;
}
}
}
GaloisFieldElement test;
for (int k_elenco=0; k_elenco<49; k_elenco++)
{
if (elenco[k_elenco].valido)
{
for (int k_E=2; k_E<14; k_E++)
{
if (E[k_E].valido)
{
test=(elenco[k_elenco].y) + (E[k_E].fuoco) * (elenco[k_elenco].x);
elenco[k_elenco].valido=found(test, E[k_E].valori);
}
}
}
}
if (conta_validi(elenco))
{
for (int i_1=0; i_1<49; i_1++)
{
if (elenco[i_1].valido)
{
if ( (elenco[i_1].x) * (elenco[i_1].y) * (a*(elenco[i_1].x) + (elenco[i_1].y)) * (c*(elenco[i_1].y) + d*(elenco[i_1].x)) * (c*(elenco[i_1].y) + e*(elenco[i_1].x)) * (f*(elenco[i_1].y) + g*(elenco[i_1].x)) * (f*(elenco[i_1].y) + h*(elenco[i_1].x)) * ((elenco[i_1].x) + (elenco[i_1].y) + uno) * (a*(elenco[i_1].x) + (elenco[i_1].x) + (elenco[i_1].y) + uno) * (c*(elenco[i_1].y) + c + d*(elenco[i_1].x) + (elenco[i_1].x)) * (c*(elenco[i_1].y) + c + e*(elenco[i_1].x) + (elenco[i_1].x)) * (f*(elenco[i_1].y) + f + g*(elenco[i_1].x) + (elenco[i_1].x)) * (f*(elenco[i_1].y) + f + h*(elenco[i_1].x) + (elenco[i_1].x)) * (a*(elenco[i_1].x) + a) * (c*(elenco[i_1].y) + d*(elenco[i_1].x) + d + (elenco[i_1].y)) * (c*(elenco[i_1].y) + e*(elenco[i_1].x) + e + (elenco[i_1].y)) * (f*(elenco[i_1].y) + g*(elenco[i_1].x) + g + (elenco[i_1].y)) * (f*(elenco[i_1].y) + h*(elenco[i_1].x) + h + (elenco[i_1].y)) * (a*c + a*(elenco[i_1].x) + c*(elenco[i_1].y) + d*(elenco[i_1].x) + d + (elenco[i_1].y)) * (a*c + a*(elenco[i_1].x) + c*(elenco[i_1].y) + e*(elenco[i_1].x) + e + (elenco[i_1].y)) * (a*f + a*(elenco[i_1].x) + f*(elenco[i_1].y) + g*(elenco[i_1].x) + g + (elenco[i_1].y)) * (a*f + a*(elenco[i_1].x) + f*(elenco[i_1].y) + h*(elenco[i_1].x) + h + (elenco[i_1].y)) * (c*d + c*e + d*(elenco[i_1].x) + e*(elenco[i_1].x)) * (c*g + c*(elenco[i_1].y) + d*f + d*(elenco[i_1].x) + f*(elenco[i_1].y) + g*(elenco[i_1].x)) * (c*h + c*(elenco[i_1].y) + d*f + d*(elenco[i_1].x) + f*(elenco[i_1].y) + h*(elenco[i_1].x)) * (c*g + c*(elenco[i_1].y) + e*f + e*(elenco[i_1].x) + f*(elenco[i_1].y) + g*(elenco[i_1].x)) * (c*h + c*(elenco[i_1].y) + e*f + e*(elenco[i_1].x) + f*(elenco[i_1].y) + h*(elenco[i_1].x)) * (f*g + f*h + g*(elenco[i_1].x) + h*(elenco[i_1].x)) != zero )
{
for (int i_2=(i_1+1); i_2<49; i_2++)
{
if (elenco[i_2].valido)
{
if ( (elenco[i_2].x) * ((elenco[i_2].y)) * (a*(elenco[i_2].x) + (elenco[i_2].y)) * (c*(elenco[i_2].y) + d*(elenco[i_2].x)) * (c*(elenco[i_2].y) + e*(elenco[i_2].x)) * (f*(elenco[i_2].y) + g*(elenco[i_2].x)) * (f*(elenco[i_2].y) + h*(elenco[i_2].x)) * ((elenco[i_1].x)*(elenco[i_2].y) + (elenco[i_1].y)*(elenco[i_2].x)) * ((elenco[i_2].x) + (elenco[i_2].y) + uno) * (a*(elenco[i_2].x) + (elenco[i_2].x) + (elenco[i_2].y) + uno) * (c*(elenco[i_2].y) + c + d*(elenco[i_2].x) + (elenco[i_2].x)) * (c*(elenco[i_2].y) + c + e*(elenco[i_2].x) + (elenco[i_2].x)) * (f*(elenco[i_2].y) + f + g*(elenco[i_2].x) + (elenco[i_2].x)) * (f*(elenco[i_2].y) + f + h*(elenco[i_2].x) + (elenco[i_2].x)) * ((elenco[i_1].x)*(elenco[i_2].y) + (elenco[i_1].x) + 

(elenco[i_1].y)*(elenco[i_2].x) + (elenco[i_2].x)) * (a*(elenco[i_2].x) + a) * (c*(elenco[i_2].y) + d*(elenco[i_2].x) + d + (elenco[i_2].y)) * (c*(elenco[i_2].y) + e*(elenco[i_2].x) + e + (elenco[i_2].y)) * (f*(elenco[i_2].y) + g*(elenco[i_2].x) + g + (elenco[i_2].y)) * (f*(elenco[i_2].y) + h*(elenco[i_2].x) + h + (elenco[i_2].y)) * ((elenco[i_1].x)*(elenco[i_2].y) + (elenco[i_1].y)*(elenco[i_2].x) + (elenco[i_1].y) + (elenco[i_2].y)) * (a*c + a*(elenco[i_2].x) + c*(elenco[i_2].y) + d*(elenco[i_2].x) + d + (elenco[i_2].y)) * (a*c + a*(elenco[i_2].x) + c*(elenco[i_2].y) + e*(elenco[i_2].x) + e + (elenco[i_2].y)) * (a*f + a*(elenco[i_2].x) + f*(elenco[i_2].y) + g*(elenco[i_2].x) + g + (elenco[i_2].y)) * (a*f + a*(elenco[i_2].x) + f*(elenco[i_2].y) + h*(elenco[i_2].x) + h + (elenco[i_2].y)) * (a*(elenco[i_1].x) + a*(elenco[i_2].x) + (elenco[i_1].x)*(elenco[i_2].y) + (elenco[i_1].y)*(elenco[i_2].x) + (elenco[i_1].y) + (elenco[i_2].y)) * (c*d + c*e + d*(elenco[i_2].x) + e*(elenco[i_2].x)) * (c*g + c*(elenco[i_2].y) + d*f + d*(elenco[i_2].x) + f*(elenco[i_2].y) + g*(elenco[i_2].x)) * (c*h + c*(elenco[i_2].y) + d*f + d*(elenco[i_2].x) + f*(elenco[i_2].y) + h*(elenco[i_2].x)) * (c*(elenco[i_1].y) + c*(elenco[i_2].y) + d*(elenco[i_1].x) + d*(elenco[i_2].x) + (elenco[i_1].x)*(elenco[i_2].y) + (elenco[i_1].y)*(elenco[i_2].x)) * (c*g + c*(elenco[i_2].y) + e*f + e*(elenco[i_2].x) + f*(elenco[i_2].y) + g*(elenco[i_2].x)) * (c*h + c*(elenco[i_2].y) + e*f + e*(elenco[i_2].x) + f*(elenco[i_2].y) + h*(elenco[i_2].x)) * (c*(elenco[i_1].y) + c*(elenco[i_2].y) + e*(elenco[i_1].x) + e*(elenco[i_2].x) + (elenco[i_1].x)*(elenco[i_2].y) + (elenco[i_1].y)*(elenco[i_2].x)) * (f*g + f*h + g*(elenco[i_2].x) + h*(elenco[i_2].x)) * (f*(elenco[i_1].y) + f*(elenco[i_2].y) + g*(elenco[i_1].x) + g*(elenco[i_2].x) + (elenco[i_1].x)*(elenco[i_2].y) + (elenco[i_1].y)*(elenco[i_2].x)) * (f*(elenco[i_1].y) + f*(elenco[i_2].y) + h*(elenco[i_1].x) + h*(elenco[i_2].x) + (elenco[i_1].x)*(elenco[i_2].y) + (elenco[i_1].y)*(elenco[i_2].x)) != zero )
{
for (int i_3=(i_2+1); i_3<49; i_3++)
{
if (elenco[i_3].valido)
{
if ( (elenco[i_3].x) * (elenco[i_3].y) * (a*(elenco[i_3].x) + (elenco[i_3].y)) * (c*(elenco[i_3].y) + d*(elenco[i_3].x)) * (c*(elenco[i_3].y) + e*(elenco[i_3].x)) * (f*(elenco[i_3].y) + g*(elenco[i_3].x)) * (f*(elenco[i_3].y) + h*(elenco[i_3].x)) * ((elenco[i_1].x)*(elenco[i_3].y) + (elenco[i_1].y)*(elenco[i_3].x)) * ((elenco[i_2].x)*(elenco[i_3].y) + (elenco[i_2].y)*(elenco[i_3].x)) * ((elenco[i_3].x) + (elenco[i_3].y) + uno) * (a*(elenco[i_3].x) + (elenco[i_3].x) + (elenco[i_3].y) + uno) * (c*(elenco[i_3].y) + c + d*(elenco[i_3].x) + (elenco[i_3].x)) * (c*(elenco[i_3].y) + c + e*(elenco[i_3].x) + (elenco[i_3].x)) * (f*(elenco[i_3].y) + f + g*(elenco[i_3].x) + (elenco[i_3].x)) * (f*(elenco[i_3].y) + f + h*(elenco[i_3].x) + (elenco[i_3].x)) * ((elenco[i_1].x)*(elenco[i_3].y) + (elenco[i_1].x) + (elenco[i_1].y)*(elenco[i_3].x) + (elenco[i_3].x)) * ((elenco[i_2].x)*(elenco[i_3].y) + (elenco[i_2].x) + (elenco[i_2].y)*(elenco[i_3].x) + (elenco[i_3].x)) * (a*(elenco[i_3].x) + a) * (c*(elenco[i_3].y) + d*(elenco[i_3].x) + d + (elenco[i_3].y)) * (c*(elenco[i_3].y) + e*(elenco[i_3].x) + e + (elenco[i_3].y)) * (f*(elenco[i_3].y) + g*(elenco[i_3].x) + g + (elenco[i_3].y)) * (f*(elenco[i_3].y) + h*(elenco[i_3].x) + h + (elenco[i_3].y)) * ((elenco[i_1].x)*(elenco[i_3].y) + (elenco[i_1].y)*(elenco[i_3].x) + (elenco[i_1].y) + (elenco[i_3].y)) * ((elenco[i_2].x)*(elenco[i_3].y) + (elenco[i_2].y)*(elenco[i_3].x) + (elenco[i_2].y) + (elenco[i_3].y)) * (a*c + a*(elenco[i_3].x) + c*(elenco[i_3].y) + d*(elenco[i_3].x) + d + 

(elenco[i_3].y)) * (a*c + a*(elenco[i_3].x) + c*(elenco[i_3].y) + e*(elenco[i_3].x) + e + (elenco[i_3].y)) * (a*f + a*(elenco[i_3].x) + f*(elenco[i_3].y) + g*(elenco[i_3].x) + g + (elenco[i_3].y)) * (a*f + a*(elenco[i_3].x) + f*(elenco[i_3].y) + h*(elenco[i_3].x) + h + (elenco[i_3].y)) * (a*(elenco[i_1].x) + a*(elenco[i_3].x) + (elenco[i_1].x)*(elenco[i_3].y) + (elenco[i_1].y)*(elenco[i_3].x) + (elenco[i_1].y) + (elenco[i_3].y)) * (a*(elenco[i_2].x) + a*(elenco[i_3].x) + (elenco[i_2].x)*(elenco[i_3].y) + (elenco[i_2].y)*(elenco[i_3].x) + (elenco[i_2].y) + (elenco[i_3].y)) * (c*d + c*e + d*(elenco[i_3].x) + e*(elenco[i_3].x)) * (c*g + c*(elenco[i_3].y) + d*f + d*(elenco[i_3].x) + f*(elenco[i_3].y) + g*(elenco[i_3].x)) * (c*h + c*(elenco[i_3].y) + d*f + d*(elenco[i_3].x) + f*(elenco[i_3].y) + h*(elenco[i_3].x)) * (c*(elenco[i_1].y) + c*(elenco[i_3].y) + d*(elenco[i_1].x) + d*(elenco[i_3].x) + (elenco[i_1].x)*(elenco[i_3].y) + (elenco[i_1].y)*(elenco[i_3].x)) * (c*(elenco[i_2].y) + c*(elenco[i_3].y) + d*(elenco[i_2].x) + d*(elenco[i_3].x) + (elenco[i_2].x)*(elenco[i_3].y) + (elenco[i_2].y)*(elenco[i_3].x)) * (c*g + c*(elenco[i_3].y) + e*f + e*(elenco[i_3].x) + f*(elenco[i_3].y) + g*(elenco[i_3].x)) * (c*h + c*(elenco[i_3].y) + e*f + e*(elenco[i_3].x) + f*(elenco[i_3].y) + h*(elenco[i_3].x)) * (c*(elenco[i_1].y) + c*(elenco[i_3].y) + e*(elenco[i_1].x) + e*(elenco[i_3].x) + (elenco[i_1].x)*(elenco[i_3].y) + (elenco[i_1].y)*(elenco[i_3].x)) * (c*(elenco[i_2].y) + c*(elenco[i_3].y) + e*(elenco[i_2].x) + e*(elenco[i_3].x) + (elenco[i_2].x)*(elenco[i_3].y) + (elenco[i_2].y)*(elenco[i_3].x)) * (f*g + f*h + g*(elenco[i_3].x) + h*(elenco[i_3].x)) * (f*(elenco[i_1].y) + f*(elenco[i_3].y) + g*(elenco[i_1].x) + g*(elenco[i_3].x) + (elenco[i_1].x)*(elenco[i_3].y) + (elenco[i_1].y)*(elenco[i_3].x)) * (f*(elenco[i_2].y) + f*(elenco[i_3].y) + g*(elenco[i_2].x) + g*(elenco[i_3].x) + (elenco[i_2].x)*(elenco[i_3].y) + (elenco[i_2].y)*(elenco[i_3].x)) * (f*(elenco[i_1].y) + f*(elenco[i_3].y) + h*(elenco[i_1].x) + h*(elenco[i_3].x) + (elenco[i_1].x)*(elenco[i_3].y) + (elenco[i_1].y)*(elenco[i_3].x)) * (f*(elenco[i_2].y) + f*(elenco[i_3].y) + h*(elenco[i_2].x) + h*(elenco[i_3].x) + (elenco[i_2].x)*(elenco[i_3].y) + (elenco[i_2].y)*(elenco[i_3].x)) * ((elenco[i_1].x)*(elenco[i_2].y) + (elenco[i_1].x)*(elenco[i_3].y) + (elenco[i_1].y)*(elenco[i_2].x) + (elenco[i_1].y)*(elenco[i_3].x) + (elenco[i_2].x)*(elenco[i_3].y) + (elenco[i_2].y)*(elenco[i_3].x)) != zero )
{
for (int i_4=(i_3+1);i_4<49;i_4++)
{
if (elenco[i_4].valido)
{
if ( (elenco[i_4].x) * (elenco[i_4].y) * (a*(elenco[i_4].x) + (elenco[i_4].y)) * (c*(elenco[i_4].y) + d*(elenco[i_4].x)) * (c*(elenco[i_4].y) + e*(elenco[i_4].x)) * (f*(elenco[i_4].y) + g*(elenco[i_4].x)) * (f*(elenco[i_4].y) + h*(elenco[i_4].x)) * ((elenco[i_1].x)*(elenco[i_4].y) + (elenco[i_1].y)*(elenco[i_4].x)) * ((elenco[i_2].x)*(elenco[i_4].y) + (elenco[i_2].y)*(elenco[i_4].x)) * ((elenco[i_3].x)*(elenco[i_4].y) + (elenco[i_3].y)*(elenco[i_4].x)) * ((elenco[i_4].x) + (elenco[i_4].y) + uno) * (a*(elenco[i_4].x) + (elenco[i_4].x) + (elenco[i_4].y) + uno) * (c*(elenco[i_4].y) + c + d*(elenco[i_4].x) + (elenco[i_4].x)) * (c*(elenco[i_4].y) + c + e*(elenco[i_4].x) + (elenco[i_4].x)) * (f*(elenco[i_4].y) + f + g*(elenco[i_4].x) + (elenco[i_4].x)) * (f*(elenco[i_4].y) + f + h*(elenco[i_4].x) + (elenco[i_4].x)) * ((elenco[i_1].x)*(elenco[i_4].y) + (elenco[i_1].x) + (elenco[i_1].y)*(elenco[i_4].x) + (elenco[i_4].x)) * ((elenco[i_2].x)*(elenco[i_4].y) + (elenco[i_2].x) + 

(elenco[i_2].y)*(elenco[i_4].x) + (elenco[i_4].x)) * ((elenco[i_3].x)*(elenco[i_4].y) + (elenco[i_3].x) + (elenco[i_3].y)*(elenco[i_4].x) + (elenco[i_4].x)) * (a*(elenco[i_4].x) + a) * (c*(elenco[i_4].y) + d*(elenco[i_4].x) + d + (elenco[i_4].y)) * (c*(elenco[i_4].y) + e*(elenco[i_4].x) + e + (elenco[i_4].y)) * (f*(elenco[i_4].y) + g*(elenco[i_4].x) + g + (elenco[i_4].y)) * (f*(elenco[i_4].y) + h*(elenco[i_4].x) + h + (elenco[i_4].y)) * ((elenco[i_1].x)*(elenco[i_4].y) + (elenco[i_1].y)*(elenco[i_4].x) + (elenco[i_1].y) + (elenco[i_4].y)) * ((elenco[i_2].x)*(elenco[i_4].y) + (elenco[i_2].y)*(elenco[i_4].x) + (elenco[i_2].y) + (elenco[i_4].y)) * ((elenco[i_3].x)*(elenco[i_4].y) + (elenco[i_3].y)*(elenco[i_4].x) + (elenco[i_3].y) + (elenco[i_4].y)) * (a*c + a*(elenco[i_4].x) + c*(elenco[i_4].y) + d*(elenco[i_4].x) + d + (elenco[i_4].y)) * (a*c + a*(elenco[i_4].x) + c*(elenco[i_4].y) + e*(elenco[i_4].x) + e + (elenco[i_4].y)) * (a*f + a*(elenco[i_4].x) + f*(elenco[i_4].y) + g*(elenco[i_4].x) + g + (elenco[i_4].y)) * (a*f + a*(elenco[i_4].x) + f*(elenco[i_4].y) + h*(elenco[i_4].x) + h + (elenco[i_4].y)) * (a*(elenco[i_1].x) + a*(elenco[i_4].x) + (elenco[i_1].x)*(elenco[i_4].y) + (elenco[i_1].y)*(elenco[i_4].x) + (elenco[i_1].y) + (elenco[i_4].y)) * (a*(elenco[i_2].x) + a*(elenco[i_4].x) + (elenco[i_2].x)*(elenco[i_4].y) + (elenco[i_2].y)*(elenco[i_4].x) + (elenco[i_2].y) + (elenco[i_4].y)) * (a*(elenco[i_3].x) + a*(elenco[i_4].x) + (elenco[i_3].x)*(elenco[i_4].y) + (elenco[i_3].y)*(elenco[i_4].x) + (elenco[i_3].y) + (elenco[i_4].y)) * (c*d + c*e + d*(elenco[i_4].x) + e*(elenco[i_4].x)) * (c*g + c*(elenco[i_4].y) + d*f + d*(elenco[i_4].x) + f*(elenco[i_4].y) + g*(elenco[i_4].x)) * (c*h + c*(elenco[i_4].y) + d*f + d*(elenco[i_4].x) + f*(elenco[i_4].y) + h*(elenco[i_4].x)) * (c*(elenco[i_1].y) + c*(elenco[i_4].y) + d*(elenco[i_1].x) + d*(elenco[i_4].x) + (elenco[i_1].x)*(elenco[i_4].y) + (elenco[i_1].y)*(elenco[i_4].x)) * (c*(elenco[i_2].y) + c*(elenco[i_4].y) + d*(elenco[i_2].x) + d*(elenco[i_4].x) + (elenco[i_2].x)*(elenco[i_4].y) + (elenco[i_2].y)*(elenco[i_4].x)) * (c*(elenco[i_3].y) + c*(elenco[i_4].y) + d*(elenco[i_3].x) + d*(elenco[i_4].x) + (elenco[i_3].x)*(elenco[i_4].y) + (elenco[i_3].y)*(elenco[i_4].x)) * (c*g + c*(elenco[i_4].y) + e*f + e*(elenco[i_4].x) + f*(elenco[i_4].y) + g*(elenco[i_4].x)) * (c*h + c*(elenco[i_4].y) + e*f + e*(elenco[i_4].x) + f*(elenco[i_4].y) + h*(elenco[i_4].x)) * (c*(elenco[i_1].y) + c*(elenco[i_4].y) + e*(elenco[i_1].x) + e*(elenco[i_4].x) + (elenco[i_1].x)*(elenco[i_4].y) + 

(elenco[i_1].y)*(elenco[i_4].x)) * (c*(elenco[i_2].y) + c*(elenco[i_4].y) + e*(elenco[i_2].x) + e*(elenco[i_4].x) + (elenco[i_2].x)*(elenco[i_4].y) + (elenco[i_2].y)*(elenco[i_4].x)) * (c*(elenco[i_3].y) + c*(elenco[i_4].y) + e*(elenco[i_3].x) + e*(elenco[i_4].x) + (elenco[i_3].x)*(elenco[i_4].y) + (elenco[i_3].y)*(elenco[i_4].x)) * (f*g + f*h + g*(elenco[i_4].x) + h*(elenco[i_4].x)) * (f*(elenco[i_1].y) + f*(elenco[i_4].y) + g*(elenco[i_1].x) + g*(elenco[i_4].x) + (elenco[i_1].x)*(elenco[i_4].y) + (elenco[i_1].y)*(elenco[i_4].x)) * (f*(elenco[i_2].y) + f*(elenco[i_4].y) + g*(elenco[i_2].x) + g*(elenco[i_4].x) + (elenco[i_2].x)*(elenco[i_4].y) + (elenco[i_2].y)*(elenco[i_4].x)) * (f*(elenco[i_3].y) + f*(elenco[i_4].y) + g*(elenco[i_3].x) + g*(elenco[i_4].x) + (elenco[i_3].x)*(elenco[i_4].y) + (elenco[i_3].y)*(elenco[i_4].x)) * (f*(elenco[i_1].y) + f*(elenco[i_4].y) + h*(elenco[i_1].x) + h*(elenco[i_4].x) + (elenco[i_1].x)*(elenco[i_4].y) + (elenco[i_1].y)*(elenco[i_4].x)) * (f*(elenco[i_2].y) + f*(elenco[i_4].y) + h*(elenco[i_2].x) + h*(elenco[i_4].x) + (elenco[i_2].x)*(elenco[i_4].y) + (elenco[i_2].y)*(elenco[i_4].x)) * (f*(elenco[i_3].y) + f*(elenco[i_4].y) + h*(elenco[i_3].x) + h*(elenco[i_4].x) + (elenco[i_3].x)*(elenco[i_4].y) + (elenco[i_3].y)*(elenco[i_4].x)) * ((elenco[i_1].x)*(elenco[i_2].y) + (elenco[i_1].x)*(elenco[i_4].y) + (elenco[i_1].y)*(elenco[i_2].x) + (elenco[i_1].y)*(elenco[i_4].x) + (elenco[i_2].x)*(elenco[i_4].y) + (elenco[i_2].y)*(elenco[i_4].x)) * ((elenco[i_1].x)*(elenco[i_3].y) + (elenco[i_1].x)*(elenco[i_4].y) + (elenco[i_1].y)*(elenco[i_3].x) + (elenco[i_1].y)*(elenco[i_4].x) + (elenco[i_3].x)*(elenco[i_4].y) + (elenco[i_3].y)*(elenco[i_4].x)) * ((elenco[i_2].x)*(elenco[i_3].y) + (elenco[i_2].x)*(elenco[i_4].y) + (elenco[i_2].y)*(elenco[i_3].x) + (elenco[i_2].y)*(elenco[i_4].x) + (elenco[i_3].x)*(elenco[i_4].y) + (elenco[i_3].y)*(elenco[i_4].x)) != zero )
{
for (int i_5=(i_4+1);i_5<49;i_5++)
{
if (elenco[i_5].valido)
{
if ( (elenco[i_5].x) * (elenco[i_5].y) * (a*(elenco[i_5].x) + (elenco[i_5].y)) * (c*(elenco[i_5].y) + d*(elenco[i_5].x)) * (c*(elenco[i_5].y) + e*(elenco[i_5].x)) * (f*(elenco[i_5].y) + g*(elenco[i_5].x)) * (f*(elenco[i_5].y) + h*(elenco[i_5].x)) * ((elenco[i_1].x)*(elenco[i_5].y) + (elenco[i_1].y)*(elenco[i_5].x)) * ((elenco[i_2].x)*(elenco[i_5].y) + (elenco[i_2].y)*(elenco[i_5].x)) * ((elenco[i_3].x)*(elenco[i_5].y) + (elenco[i_3].y)*(elenco[i_5].x)) * ((elenco[i_4].x)*(elenco[i_5].y) + (elenco[i_4].y)*(elenco[i_5].x)) * ((elenco[i_5].x) + (elenco[i_5].y) + uno) * (a*(elenco[i_5].x) + (elenco[i_5].x) + (elenco[i_5].y) + uno) * (c*(elenco[i_5].y) + c + d*(elenco[i_5].x) + (elenco[i_5].x)) * (c*(elenco[i_5].y) + c + e*(elenco[i_5].x) + (elenco[i_5].x)) * (f*(elenco[i_5].y) + f + g*(elenco[i_5].x) + (elenco[i_5].x)) * (f*(elenco[i_5].y) + f + h*(elenco[i_5].x) + (elenco[i_5].x)) * ((elenco[i_1].x)*(elenco[i_5].y) + (elenco[i_1].x) + (elenco[i_1].y)*(elenco[i_5].x) + (elenco[i_5].x)) * ((elenco[i_2].x)*(elenco[i_5].y) + (elenco[i_2].x) + (elenco[i_2].y)*(elenco[i_5].x) + (elenco[i_5].x)) * ((elenco[i_3].x)*(elenco[i_5].y) + (elenco[i_3].x) + (elenco[i_3].y)*(elenco[i_5].x) + (elenco[i_5].x)) * ((elenco[i_4].x)*(elenco[i_5].y) + (elenco[i_4].x) + (elenco[i_4].y)*(elenco[i_5].x) + (elenco[i_5].x)) * (a*(elenco[i_5].x) + a) * (c*(elenco[i_5].y) + d*(elenco[i_5].x) + d + (elenco[i_5].y)) * (c*(elenco[i_5].y) + e*(elenco[i_5].x) + e + (elenco[i_5].y)) * (f*(elenco[i_5].y) + g*(elenco[i_5].x) + g + (elenco[i_5].y)) * (f*(elenco[i_5].y) + h*(elenco[i_5].x) + h + (elenco[i_5].y)) * ((elenco[i_1].x)*(elenco[i_5].y) + (elenco[i_1].y)*(elenco[i_5].x) + (elenco[i_1].y) + (elenco[i_5].y)) * ((elenco[i_2].x)*(elenco[i_5].y) + (elenco[i_2].y)*(elenco[i_5].x) + (elenco[i_2].y) + (elenco[i_5].y)) * ((elenco[i_3].x)*(elenco[i_5].y) + (elenco[i_3].y)*(elenco[i_5].x) + (elenco[i_3].y) + (elenco[i_5].y)) * ((elenco[i_4].x)*(elenco[i_5].y) + (elenco[i_4].y)*(elenco[i_5].x) + (elenco[i_4].y) + (elenco[i_5].y)) * (a*c + a*(elenco[i_5].x) + c*(elenco[i_5].y) + d*(elenco[i_5].x) + d + (elenco[i_5].y)) * (a*c + a*(elenco[i_5].x) + c*(elenco[i_5].y) + e*(elenco[i_5].x) + e + (elenco[i_5].y)) * (a*f + a*(elenco[i_5].x) + f*(elenco[i_5].y) + g*(elenco[i_5].x) + g + (elenco[i_5].y)) * (a*f + a*(elenco[i_5].x) + f*(elenco[i_5].y) + h*(elenco[i_5].x) + h + (elenco[i_5].y)) * (a*(elenco[i_1].x) + a*(elenco[i_5].x) + (elenco[i_1].x)*(elenco[i_5].y) + (elenco[i_1].y)*(elenco[i_5].x) + 

(elenco[i_1].y) + (elenco[i_5].y)) * (a*(elenco[i_2].x) + a*(elenco[i_5].x) + (elenco[i_2].x)*(elenco[i_5].y) + (elenco[i_2].y)*(elenco[i_5].x) + (elenco[i_2].y) + (elenco[i_5].y)) * (a*(elenco[i_3].x) + a*(elenco[i_5].x) + (elenco[i_3].x)*(elenco[i_5].y) + (elenco[i_3].y)*(elenco[i_5].x) + (elenco[i_3].y) + (elenco[i_5].y)) * (a*(elenco[i_4].x) + a*(elenco[i_5].x) + (elenco[i_4].x)*(elenco[i_5].y) + (elenco[i_4].y)*(elenco[i_5].x) + (elenco[i_4].y) + (elenco[i_5].y)) * (c*d + c*e + d*(elenco[i_5].x) + e*(elenco[i_5].x)) * (c*g + c*(elenco[i_5].y) + d*f + d*(elenco[i_5].x) + f*(elenco[i_5].y) + g*(elenco[i_5].x)) * (c*h + c*(elenco[i_5].y) + d*f + d*(elenco[i_5].x) + f*(elenco[i_5].y) + h*(elenco[i_5].x)) * (c*(elenco[i_1].y) + c*(elenco[i_5].y) + d*(elenco[i_1].x) + d*(elenco[i_5].x) + (elenco[i_1].x)*(elenco[i_5].y) + (elenco[i_1].y)*(elenco[i_5].x)) * (c*(elenco[i_2].y) + c*(elenco[i_5].y) + d*(elenco[i_2].x) + d*(elenco[i_5].x) + (elenco[i_2].x)*(elenco[i_5].y) + (elenco[i_2].y)*(elenco[i_5].x)) * (c*(elenco[i_3].y) + c*(elenco[i_5].y) + d*(elenco[i_3].x) + d*(elenco[i_5].x) + (elenco[i_3].x)*(elenco[i_5].y) + (elenco[i_3].y)*(elenco[i_5].x)) * (c*(elenco[i_4].y) + c*(elenco[i_5].y) + d*(elenco[i_4].x) + d*(elenco[i_5].x) + (elenco[i_4].x)*(elenco[i_5].y) + (elenco[i_4].y)*(elenco[i_5].x)) * (c*g + c*(elenco[i_5].y) + e*f + e*(elenco[i_5].x) + f*(elenco[i_5].y) + g*(elenco[i_5].x)) * (c*h + c*(elenco[i_5].y) + e*f + e*(elenco[i_5].x) + f*(elenco[i_5].y) + h*(elenco[i_5].x)) * (c*(elenco[i_1].y) + c*(elenco[i_5].y) + e*(elenco[i_1].x) + e*(elenco[i_5].x) + (elenco[i_1].x)*(elenco[i_5].y) + (elenco[i_1].y)*(elenco[i_5].x)) * (c*(elenco[i_2].y) + c*(elenco[i_5].y) + e*(elenco[i_2].x) + e*(elenco[i_5].x) + (elenco[i_2].x)*(elenco[i_5].y) + (elenco[i_2].y)*(elenco[i_5].x)) * (c*(elenco[i_3].y) + c*(elenco[i_5].y) + e*(elenco[i_3].x) + e*(elenco[i_5].x) + (elenco[i_3].x)*(elenco[i_5].y) + (elenco[i_3].y)*(elenco[i_5].x)) * (c*(elenco[i_4].y) + c*(elenco[i_5].y) + e*(elenco[i_4].x) + e*(elenco[i_5].x) + (elenco[i_4].x)*(elenco[i_5].y) + (elenco[i_4].y)*(elenco[i_5].x)) * (f*g + f*h + g*(elenco[i_5].x) + h*(elenco[i_5].x)) * (f*(elenco[i_1].y) + f*(elenco[i_5].y) + g*(elenco[i_1].x) + g*(elenco[i_5].x) + (elenco[i_1].x)*(elenco[i_5].y) + 

(elenco[i_1].y)*(elenco[i_5].x)) * (f*(elenco[i_2].y) + f*(elenco[i_5].y) + g*(elenco[i_2].x) + g*(elenco[i_5].x) + (elenco[i_2].x)*(elenco[i_5].y) + (elenco[i_2].y)*(elenco[i_5].x)) * (f*(elenco[i_3].y) + f*(elenco[i_5].y) + g*(elenco[i_3].x) + g*(elenco[i_5].x) + (elenco[i_3].x)*(elenco[i_5].y) + (elenco[i_3].y)*(elenco[i_5].x)) * (f*(elenco[i_4].y) + f*(elenco[i_5].y) + g*(elenco[i_4].x) + g*(elenco[i_5].x) + (elenco[i_4].x)*(elenco[i_5].y) + (elenco[i_4].y)*(elenco[i_5].x)) * (f*(elenco[i_1].y) + f*(elenco[i_5].y) + h*(elenco[i_1].x) + h*(elenco[i_5].x) + (elenco[i_1].x)*(elenco[i_5].y) + (elenco[i_1].y)*(elenco[i_5].x)) * (f*(elenco[i_2].y) + f*(elenco[i_5].y) + h*(elenco[i_2].x) + h*(elenco[i_5].x) + (elenco[i_2].x)*(elenco[i_5].y) + (elenco[i_2].y)*(elenco[i_5].x)) * (f*(elenco[i_3].y) + f*(elenco[i_5].y) + h*(elenco[i_3].x) + h*(elenco[i_5].x) + (elenco[i_3].x)*(elenco[i_5].y) + (elenco[i_3].y)*(elenco[i_5].x)) * (f*(elenco[i_4].y) + f*(elenco[i_5].y) + h*(elenco[i_4].x) + h*(elenco[i_5].x) + (elenco[i_4].x)*(elenco[i_5].y) + (elenco[i_4].y)*(elenco[i_5].x)) * ((elenco[i_1].x)*(elenco[i_2].y) + (elenco[i_1].x)*(elenco[i_5].y) + (elenco[i_1].y)*(elenco[i_2].x) + (elenco[i_1].y)*(elenco[i_5].x) + (elenco[i_2].x)*(elenco[i_5].y) + (elenco[i_2].y)*(elenco[i_5].x)) * ((elenco[i_1].x)*(elenco[i_3].y) + (elenco[i_1].x)*(elenco[i_5].y) + (elenco[i_1].y)*(elenco[i_3].x) + (elenco[i_1].y)*(elenco[i_5].x) + (elenco[i_3].x)*(elenco[i_5].y) + (elenco[i_3].y)*(elenco[i_5].x)) * ((elenco[i_1].x)*(elenco[i_4].y) + (elenco[i_1].x)*(elenco[i_5].y) + (elenco[i_1].y)*(elenco[i_4].x) + (elenco[i_1].y)*(elenco[i_5].x) + (elenco[i_4].x)*(elenco[i_5].y) + (elenco[i_4].y)*(elenco[i_5].x)) * ((elenco[i_2].x)*(elenco[i_3].y) + (elenco[i_2].x)*(elenco[i_5].y) + (elenco[i_2].y)*(elenco[i_3].x) + (elenco[i_2].y)*(elenco[i_5].x) + (elenco[i_3].x)*(elenco[i_5].y) + (elenco[i_3].y)*(elenco[i_5].x)) * ((elenco[i_2].x)*(elenco[i_4].y) + (elenco[i_2].x)*(elenco[i_5].y) + (elenco[i_2].y)*(elenco[i_4].x) + (elenco[i_2].y)*(elenco[i_5].x) + (elenco[i_4].x)*(elenco[i_5].y) + (elenco[i_4].y)*(elenco[i_5].x)) * ((elenco[i_3].x)*(elenco[i_4].y) + (elenco[i_3].x)*(elenco[i_5].y) + (elenco[i_3].y)*(elenco[i_4].x) + (elenco[i_3].y)*(elenco[i_5].x) + (elenco[i_4].x)*(elenco[i_5].y) + (elenco[i_4].y)*(elenco[i_5].x)) != zero )
{
for (int i_6=(i_5+1);i_6<49; i_6++)
{
if (elenco[i_6].valido)
{
if ( (elenco[i_6].x) * (elenco[i_6].y) * (a*(elenco[i_6].x) + (elenco[i_6].y)) * (c*(elenco[i_6].y) + d*(elenco[i_6].x)) * (c*(elenco[i_6].y) + e*(elenco[i_6].x)) * (f*(elenco[i_6].y) + g*(elenco[i_6].x)) * (f*(elenco[i_6].y) + h*(elenco[i_6].x)) * ((elenco[i_1].x)*(elenco[i_6].y) + (elenco[i_1].y)*(elenco[i_6].x)) * ((elenco[i_2].x)*(elenco[i_6].y) + (elenco[i_2].y)*(elenco[i_6].x)) * ((elenco[i_3].x)*(elenco[i_6].y) + (elenco[i_3].y)*(elenco[i_6].x)) * ((elenco[i_4].x)*(elenco[i_6].y) + (elenco[i_4].y)*(elenco[i_6].x)) * ((elenco[i_5].x)*(elenco[i_6].y) + (elenco[i_5].y)*(elenco[i_6].x)) * ((elenco[i_6].x) + (elenco[i_6].y) + uno) * (a*(elenco[i_6].x) + (elenco[i_6].x) + (elenco[i_6].y) + uno) * (c*(elenco[i_6].y) + c + d*(elenco[i_6].x) + (elenco[i_6].x)) * (c*(elenco[i_6].y) + c + e*(elenco[i_6].x) + (elenco[i_6].x)) * (f*(elenco[i_6].y) + f + g*(elenco[i_6].x) + (elenco[i_6].x)) * (f*(elenco[i_6].y) + f + h*(elenco[i_6].x) + (elenco[i_6].x)) * ((elenco[i_1].x)*(elenco[i_6].y) + (elenco[i_1].x) + (elenco[i_1].y)*(elenco[i_6].x) + (elenco[i_6].x)) * ((elenco[i_2].x)*(elenco[i_6].y) + (elenco[i_2].x) + (elenco[i_2].y)*(elenco[i_6].x) + (elenco[i_6].x)) * ((elenco[i_3].x)*(elenco[i_6].y) + (elenco[i_3].x) + (elenco[i_3].y)*(elenco[i_6].x) + (elenco[i_6].x)) * ((elenco[i_4].x)*(elenco[i_6].y) + (elenco[i_4].x) + (elenco[i_4].y)*(elenco[i_6].x) + (elenco[i_6].x)) * ((elenco[i_5].x)*(elenco[i_6].y) + (elenco[i_5].x) + (elenco[i_5].y)*(elenco[i_6].x) + (elenco[i_6].x)) * (a*(elenco[i_6].x) + a) * (c*(elenco[i_6].y) + d*(elenco[i_6].x) + d + (elenco[i_6].y)) * (c*(elenco[i_6].y) + e*(elenco[i_6].x) + e + (elenco[i_6].y)) * (f*(elenco[i_6].y) + g*(elenco[i_6].x) + g + 

(elenco[i_6].y)) * (f*(elenco[i_6].y) + h*(elenco[i_6].x) + h + (elenco[i_6].y)) * ((elenco[i_1].x)*(elenco[i_6].y) + (elenco[i_1].y)*(elenco[i_6].x) + (elenco[i_1].y) + (elenco[i_6].y)) * ((elenco[i_2].x)*(elenco[i_6].y) + (elenco[i_2].y)*(elenco[i_6].x) + (elenco[i_2].y) + (elenco[i_6].y)) * ((elenco[i_3].x)*(elenco[i_6].y) + (elenco[i_3].y)*(elenco[i_6].x) + (elenco[i_3].y) + (elenco[i_6].y)) * ((elenco[i_4].x)*(elenco[i_6].y) + (elenco[i_4].y)*(elenco[i_6].x) + (elenco[i_4].y) + (elenco[i_6].y)) * ((elenco[i_5].x)*(elenco[i_6].y) + (elenco[i_5].y)*(elenco[i_6].x) + (elenco[i_5].y) + (elenco[i_6].y)) * (a*c + a*(elenco[i_6].x) + c*(elenco[i_6].y) + d*(elenco[i_6].x) + d + (elenco[i_6].y)) * (a*c + a*(elenco[i_6].x) + c*(elenco[i_6].y) + e*(elenco[i_6].x) + e + (elenco[i_6].y)) * (a*f + a*(elenco[i_6].x) + f*(elenco[i_6].y) + g*(elenco[i_6].x) + g + (elenco[i_6].y)) * (a*f + a*(elenco[i_6].x) + f*(elenco[i_6].y) + h*(elenco[i_6].x) + h + (elenco[i_6].y)) * (a*(elenco[i_1].x) + a*(elenco[i_6].x) + (elenco[i_1].x)*(elenco[i_6].y) + (elenco[i_1].y)*(elenco[i_6].x) + (elenco[i_1].y) + (elenco[i_6].y)) * (a*(elenco[i_2].x) + a*(elenco[i_6].x) + (elenco[i_2].x)*(elenco[i_6].y) + (elenco[i_2].y)*(elenco[i_6].x) + (elenco[i_2].y) + (elenco[i_6].y)) * (a*(elenco[i_3].x) + a*(elenco[i_6].x) + (elenco[i_3].x)*(elenco[i_6].y) + (elenco[i_3].y)*(elenco[i_6].x) + (elenco[i_3].y) + (elenco[i_6].y)) * (a*(elenco[i_4].x) + a*(elenco[i_6].x) + (elenco[i_4].x)*(elenco[i_6].y) + (elenco[i_4].y)*(elenco[i_6].x) + (elenco[i_4].y) + (elenco[i_6].y)) * (a*(elenco[i_5].x) + a*(elenco[i_6].x) + (elenco[i_5].x)*(elenco[i_6].y) + (elenco[i_5].y)*(elenco[i_6].x) + (elenco[i_5].y) + (elenco[i_6].y)) * (c*d + c*e + d*(elenco[i_6].x) + e*(elenco[i_6].x)) * (c*g + c*(elenco[i_6].y) + d*f + d*(elenco[i_6].x) + f*(elenco[i_6].y) + g*(elenco[i_6].x)) * (c*h + c*(elenco[i_6].y) + d*f + d*(elenco[i_6].x) + f*(elenco[i_6].y) + h*(elenco[i_6].x)) * (c*(elenco[i_1].y) + c*(elenco[i_6].y) + d*(elenco[i_1].x) + d*(elenco[i_6].x) + (elenco[i_1].x)*(elenco[i_6].y) + (elenco[i_1].y)*(elenco[i_6].x)) * (c*(elenco[i_2].y) + c*(elenco[i_6].y) + 

d*(elenco[i_2].x) + d*(elenco[i_6].x) + (elenco[i_2].x)*(elenco[i_6].y) + (elenco[i_2].y)*(elenco[i_6].x)) * (c*(elenco[i_3].y) + c*(elenco[i_6].y) + d*(elenco[i_3].x) + d*(elenco[i_6].x) + (elenco[i_3].x)*(elenco[i_6].y) + (elenco[i_3].y)*(elenco[i_6].x)) * (c*(elenco[i_4].y) + c*(elenco[i_6].y) + d*(elenco[i_4].x) + d*(elenco[i_6].x) + (elenco[i_4].x)*(elenco[i_6].y) + (elenco[i_4].y)*(elenco[i_6].x)) * (c*(elenco[i_5].y) + c*(elenco[i_6].y) + d*(elenco[i_5].x) + d*(elenco[i_6].x) + (elenco[i_5].x)*(elenco[i_6].y) + (elenco[i_5].y)*(elenco[i_6].x)) * (c*g + c*(elenco[i_6].y) + e*f + e*(elenco[i_6].x) + f*(elenco[i_6].y) + g*(elenco[i_6].x)) * (c*h + c*(elenco[i_6].y) + e*f + e*(elenco[i_6].x) + f*(elenco[i_6].y) + h*(elenco[i_6].x)) * (c*(elenco[i_1].y) + c*(elenco[i_6].y) + e*(elenco[i_1].x) + e*(elenco[i_6].x) + (elenco[i_1].x)*(elenco[i_6].y) + (elenco[i_1].y)*(elenco[i_6].x)) * (c*(elenco[i_2].y) + c*(elenco[i_6].y) + e*(elenco[i_2].x) + e*(elenco[i_6].x) + (elenco[i_2].x)*(elenco[i_6].y) + (elenco[i_2].y)*(elenco[i_6].x)) * (c*(elenco[i_3].y) + c*(elenco[i_6].y) + e*(elenco[i_3].x) + e*(elenco[i_6].x) + (elenco[i_3].x)*(elenco[i_6].y) + (elenco[i_3].y)*(elenco[i_6].x)) * (c*(elenco[i_4].y) + c*(elenco[i_6].y) + e*(elenco[i_4].x) + e*(elenco[i_6].x) + (elenco[i_4].x)*(elenco[i_6].y) + (elenco[i_4].y)*(elenco[i_6].x)) * (c*(elenco[i_5].y) + c*(elenco[i_6].y) + e*(elenco[i_5].x) + e*(elenco[i_6].x) + (elenco[i_5].x)*(elenco[i_6].y) + (elenco[i_5].y)*(elenco[i_6].x)) * (f*g + f*h + g*(elenco[i_6].x) + h*(elenco[i_6].x)) * (f*(elenco[i_1].y) + f*(elenco[i_6].y) + g*(elenco[i_1].x) + g*(elenco[i_6].x) + (elenco[i_1].x)*(elenco[i_6].y) + (elenco[i_1].y)*(elenco[i_6].x)) * (f*(elenco[i_2].y) + f*(elenco[i_6].y) + g*(elenco[i_2].x) + g*(elenco[i_6].x) + (elenco[i_2].x)*(elenco[i_6].y) + (elenco[i_2].y)*(elenco[i_6].x)) * (f*(elenco[i_3].y) + f*(elenco[i_6].y) + g*(elenco[i_3].x) + g*(elenco[i_6].x) + (elenco[i_3].x)*(elenco[i_6].y) + (elenco[i_3].y)*(elenco[i_6].x)) * (f*(elenco[i_4].y) + f*(elenco[i_6].y) + g*(elenco[i_4].x) + g*(elenco[i_6].x) + (elenco[i_4].x)*(elenco[i_6].y) + (elenco[i_4].y)*(elenco[i_6].x)) * (f*(elenco[i_5].y) + f*(elenco[i_6].y) + g*(elenco[i_5].x) + g*(elenco[i_6].x) + (elenco[i_5].x)*(elenco[i_6].y) + (elenco[i_5].y)*(elenco[i_6].x)) * (f*(elenco[i_1].y) + f*(elenco[i_6].y) + h*(elenco[i_1].x) + h*(elenco[i_6].x) + (elenco[i_1].x)*(elenco[i_6].y) + (elenco[i_1].y)*(elenco[i_6].x)) * (f*(elenco[i_2].y) + f*(elenco[i_6].y) + h*(elenco[i_2].x) + h*(elenco[i_6].x) + (elenco[i_2].x)*(elenco[i_6].y) + (elenco[i_2].y)*(elenco[i_6].x)) * (f*(elenco[i_3].y) + f*(elenco[i_6].y) + h*(elenco[i_3].x) + h*(elenco[i_6].x) + (elenco[i_3].x)*(elenco[i_6].y) + (elenco[i_3].y)*(elenco[i_6].x)) * (f*(elenco[i_4].y) + f*(elenco[i_6].y) + h*(elenco[i_4].x) + h*(elenco[i_6].x) + (elenco[i_4].x)*(elenco[i_6].y) + (elenco[i_4].y)*(elenco[i_6].x)) * (f*(elenco[i_5].y) + f*(elenco[i_6].y) + h*(elenco[i_5].x) + h*(elenco[i_6].x) + (elenco[i_5].x)*(elenco[i_6].y) + 

(elenco[i_5].y)*(elenco[i_6].x)) * ((elenco[i_1].x)*(elenco[i_2].y) + (elenco[i_1].x)*(elenco[i_6].y) + (elenco[i_1].y)*(elenco[i_2].x) + (elenco[i_1].y)*(elenco[i_6].x) + (elenco[i_2].x)*(elenco[i_6].y) + (elenco[i_2].y)*(elenco[i_6].x)) * ((elenco[i_1].x)*(elenco[i_3].y) + (elenco[i_1].x)*(elenco[i_6].y) + (elenco[i_1].y)*(elenco[i_3].x) + (elenco[i_1].y)*(elenco[i_6].x) + (elenco[i_3].x)*(elenco[i_6].y) + (elenco[i_3].y)*(elenco[i_6].x)) * ((elenco[i_1].x)*(elenco[i_4].y) + (elenco[i_1].x)*(elenco[i_6].y) + (elenco[i_1].y)*(elenco[i_4].x) + (elenco[i_1].y)*(elenco[i_6].x) + (elenco[i_4].x)*(elenco[i_6].y) + (elenco[i_4].y)*(elenco[i_6].x)) * ((elenco[i_1].x)*(elenco[i_5].y) + (elenco[i_1].x)*(elenco[i_6].y) + (elenco[i_1].y)*(elenco[i_5].x) + (elenco[i_1].y)*(elenco[i_6].x) + (elenco[i_5].x)*(elenco[i_6].y) + (elenco[i_5].y)*(elenco[i_6].x)) * ((elenco[i_2].x)*(elenco[i_3].y) + (elenco[i_2].x)*(elenco[i_6].y) + (elenco[i_2].y)*(elenco[i_3].x) + (elenco[i_2].y)*(elenco[i_6].x) + (elenco[i_3].x)*(elenco[i_6].y) + (elenco[i_3].y)*(elenco[i_6].x)) * ((elenco[i_2].x)*(elenco[i_4].y) + (elenco[i_2].x)*(elenco[i_6].y) + (elenco[i_2].y)*(elenco[i_4].x) + (elenco[i_2].y)*(elenco[i_6].x) + (elenco[i_4].x)*(elenco[i_6].y) + (elenco[i_4].y)*(elenco[i_6].x)) * ((elenco[i_2].x)*(elenco[i_5].y) + (elenco[i_2].x)*(elenco[i_6].y) + (elenco[i_2].y)*(elenco[i_5].x) + (elenco[i_2].y)*(elenco[i_6].x) + (elenco[i_5].x)*(elenco[i_6].y) + (elenco[i_5].y)*(elenco[i_6].x)) * ((elenco[i_3].x)*(elenco[i_4].y) + (elenco[i_3].x)*(elenco[i_6].y) + (elenco[i_3].y)*(elenco[i_4].x) + (elenco[i_3].y)*(elenco[i_6].x) + (elenco[i_4].x)*(elenco[i_6].y) + (elenco[i_4].y)*(elenco[i_6].x)) * ((elenco[i_3].x)*(elenco[i_5].y) + (elenco[i_3].x)*(elenco[i_6].y) + (elenco[i_3].y)*(elenco[i_5].x) + (elenco[i_3].y)*(elenco[i_6].x) + (elenco[i_5].x)*(elenco[i_6].y) + (elenco[i_5].y)*(elenco[i_6].x)) * ((elenco[i_4].x)*(elenco[i_5].y) + (elenco[i_4].x)*(elenco[i_6].y) + (elenco[i_4].y)*(elenco[i_5].x) + (elenco[i_4].y)*(elenco[i_6].x) + (elenco[i_5].x)*(elenco[i_6].y) + (elenco[i_5].y)*(elenco[i_6].x)) != zero )
{
if (test_iperfocalizzato(zero, uno, a, c, d, e, f, g, h, elenco[i_1].x, elenco[i_1].y, elenco[i_2].x, elenco[i_2].y, elenco[i_3].x, elenco[i_3].y, elenco[i_4].x, elenco[i_4].y, elenco[i_5].x, elenco[i_5].y, elenco[i_6].x, elenco[i_6].y, lista, n_fuochi))
print(zero, uno, a, c, d, e, f, g, h, elenco[i_1].x, elenco[i_1].y, elenco[i_2].x, elenco[i_2].y, elenco[i_3].x, elenco[i_3].y, elenco[i_4].x, elenco[i_4].y, elenco[i_5].x, elenco[i_5].y, elenco[i_6].x, elenco[i_6].y, &n_archi);
}
}
}
}
}
}
}
}
}
}
}
}
}
}
}
}
}
}
}
}
}
}
}
}
}
}
}
}
}
}
}
}
}
}
}
}
time (&end);
diff = difftime(end,start);
string text_1("\nOperazione terminata.\nSelezionati ");
string text_2(" archi in ");
string text_3(" minuti.\n");
print(text_1, text_2, text_3, n_archi, diff/60);
cout << text_1 << n_archi << text_2 << (diff/60) << text_3;
return 0;
}


\end{lstlisting}

\end{document}